\documentclass[12pt]{article}
\usepackage{amsfonts}
\usepackage{latexsym, amssymb, amsmath, amscd, amsfonts, epsfig, graphicx}
\usepackage{mathrsfs}
\usepackage{color}
\usepackage{ifpdf}
\newtheorem{theorem}{Theorem}[section]

\newtheorem{lemma}[theorem]{Lemma}

\newtheorem{conjecture}[theorem]{Conjecture}

 \numberwithin{equation}{section}

\def\qed{\nopagebreak\hfill{\rule{4pt}{7pt}}}
\def\proof{\noindent {\it{Proof.} \hskip 2pt}}

\begin{document}

\begin{center}
{\large\bf  Higher Order Log-Concavity in Euler's Difference Table}
\end{center}

\begin{center}
William Y.C. Chen$^{1}$, Cindy C.Y. Gu$^{2}$, Kevin J. Ma$^{3}$ and Larry X.W. Wang$^{4}$\\[6pt]
Center for Combinatorics, LPMC-TJKLC\\
Nankai University, Tianjin 300071, P. R. China

$^{1}${ chen@nankai.edu.cn}, $^{2}${
guchunyan@cfc.nankai.edu.cn}, \\$^{3}${ majun@cfc nankai.edu.cn},
$^{4}${ wxw@cfc.nankai.edu.cn}
\end{center}

\vspace{0.3cm} \noindent{\bf Abstract.}  Let $e_{n}^k$ be the
entries in the classical Euler's difference table. We consider the
array $d_{n}^{k}=e_n^k/k!$ for $0\leq k \leq n$, where $d_n^k$ can
be interpreted as  the number of $k$-fixed-points-permutations of
$[n]$. We show that the sequence $\{d_n^k\}_{0\leq k\leq n}$ is 
$2$-log-concave and reverse ultra log-concave for any given $n$.

\noindent {\bf Keywords:} log-concavity, 2-log-concavity,
reverse ultra log-concavity, Euler's difference table

\noindent {\bf  Classification:} 05A20; 33F10

\section{Introduction}

Euler introduced the difference table
$(e_n^k)_{0\leq k\leq n}$ defined by $e_n^n=n!$ and
\begin{equation}\label{edt}
e_n^{k-1}=e_n^k-e_{n-1}^{k-1},
\end{equation} for $1\leq k\leq n$; see \cite{Cla}.
The combinatorial interpretation of the numbers $e_n^k$
 was found by Dumont and
Randrianarivony \cite{Dum}. Clarke, Han and Zeng \cite{Cla} further gave a
combinatorial interpretation of the $q$-analogue of Euler's
difference table, and this interpretation has been extended by
Faliharimalala and Zeng \cite{fal1, fal2} to the wreath product
$C_\ell\wr S_n$ of the cyclic group with the symmetric group.

It is easily seen from the recurrence (\ref{edt}) that $k!$ divides
$e_n^k$. Thus we can define the integers $d_n^k=e_n^k/k!$.
Rakotondrajao \cite{Rak1} has shown that $d_n^k$ counts the number
of $k$-fixed-points-permutations of $[n]$, where a permutation $\pi
\in \mathfrak{S}_n$ is called $k$-fixed-points-permutation if there
are no fixed points in the last $n-k$ positions and the first $k$
elements are  in different cycles. Based on this combinatorial
interpretation, Rakotondrajao \cite{Rak2} has found bijective proofs
for the following two recurrence relations for $0\leq k \leq n-1$,
\begin{eqnarray}\label{rec1}
d_n^k&=&(n-1)d_{n-1}^k+(n-k-1)d_{n-2}^k, 
\\[5pt]
\label{rec2}
d_n^k&=&nd_{n-1}^k-d_{n-2}^{k-1}.
\end{eqnarray}
Notice that $d_k^k=1$. Recently, Eriksen, Freij and W\"{a}stlund
\cite{Eri} have generalized these formulas to fixed point
$\lambda$-colored permutations. Employing \eqref{rec1} and
\eqref{rec2}, we can easily derive the following recurrence
relation for $0\leq k\leq n-1$,
\begin{equation}\label{rec3}
d_n^k=d_{n-1}^{k-1}+(n-k)d_{n-1}^k.
\end{equation}
Using the above recurrence relations (\ref{rec1}) (\ref{rec2})
and (\ref{rec3}), we shall prove that the sequence $\{d_n^k\}_{0\leq
k\leq n}$ has higher order log-concave properties. To be more
specific, we shall show that this sequence is $2$-log-concave and
reverse ultra log-concave for any $n\geq 1$.

\section{$2$-log-concavity}

In this section, we shall show that the sequence $\{d_n^k\}_{0\leq
k\leq n}$ is $2$-log-concave for any $n\geq 1$. Recall that
 a sequence $\{a_k\}_{k\geq 0}$ of real numbers is
said to be log-concave if $a_k^2\geq a_{k+1}a_{k-1}$  for all $k\geq
1$; see Stanley \cite{stanley1989} and Brenti \cite{bre}. From the
recurrence relation (\ref{rec3}), it is easy to prove by induction
that
  the sequence $\{d_n^k\}_{0\leq
k\leq n}$ is log-concave.
\begin{theorem}
For $1\leq k \leq n$, we have
$$(d_n^k)^2\geq d_n^{k-1}d_n^{k+1},$$
that is, the sequence $\{d_n^k\}_{0\leq k\leq n}$ is log-concave.
\end{theorem}

The notion of high order log-concavity was introduced by Moll
\cite{moll}; see also, \cite{KP}. Given a sequence $\{a_k\}_{k\geq
0}$, define the operator $\mathfrak{L}$  as
$\mathfrak{L}\{a_k\}=\{b_k\}$, where
$$b_k=a_k^2-a_{k-1}a_{k+1}.$$
The log-concavity of $\{a_k\}$  becomes the positivity of
$\mathfrak{L}\{a_k\}$. If the sequence $\mathfrak{L}\{a_k\}$ is not
only positive but also log-concave, then we say that $\{a_k\}$ is
$2$-log-concave. In general, we say that $\{a_k\}$ is
$l$-log-concave if $\mathfrak{L}^l\{a_k\}$ is positive, and that
$\{a_k\}$ is infinite log-concave if $\mathfrak{L}^l \{a_k\}$ is
positive for any $l\geq 1$. From numerical evidence, we pose the
following conjecture.

\begin{conjecture}
The sequence $\{d_n^k\}_{0\leq
k\leq n}$ is infinitely log-concave.
\end{conjecture}

Recently, Br\"{a}nd\'{e}n \cite{Bra} and Cardon \cite{Cardon} have independently
 proved  that if a polynomial
has only real and nonpositive zeros, then its Taylor coefficients
form an infinite log-concave sequence.  However, this is not the case of the polynomials $\sum
d_n^kx^k$. For example, for $n=2$, the polynomial  $x^2+x+1$ does not have real roots.
Nevertheless, we shall show that the sequence
$\{d_n^k\}$ is $2$-log concave in support of the general conjecture.

\begin{theorem}\label{t2log}
The sequence $\{d_n^k\}_{0\leq k\leq n}$ is $2$-log-concave. In
other words, for $n\geq 4$ and $2\leq k\leq n-2$, we have
\begin{equation}\label{2log}
\left((d_n^k)^2-d_n^{k-1}d_n^{k+1}\right)^2-\left((d_n^{k-1})^2-
d_n^{k-2}d_n^{k}\right)\left((d_n^{k+1})^2-d_n^{k}d_n^{k+2}\right)
\geq 0.
\end{equation}
\end{theorem}

The idea to prove Theorem \ref{t2log} may be described as follows.
As the first step, we reformulate the left hand side of the above
inequality (\ref{2log}) a cubic function $f$ on
$\frac{d_{n+1}^k}{d_n^k}$ by applying the recurrence relations
(\ref{rec1}), (\ref{rec2}), (\ref{rec3}) and the recurrence relation
presented in the following Lemma \ref{rr}. Then Theorem \ref{t2log}
is equivalent to the assertion that $f\geq 0$ on the interval
\[ I=[n+\frac{n-k}{n},n+\frac{n-k}{n}+\frac{n-k}{n^2}],\] since it can be
 verified that for $n\geq 4$ and $2\leq k\leq n-2$,
\begin{equation}\label{bound}
\frac{n-k}{n}\leq \frac{d_{n+1}^k}{d_n^k}\leq
n+\frac{n-k}{n}+\frac{n-k}{n^2}.
\end{equation}
Moreover,  when $f(x)$ is considered as a continuous function on $x$,
 we will be able to show
that $f'(x)<0$ for  $x\in I$ and
\[ f\left(n+\frac{n-k}{n}+\frac{n-k}{n^2}\right)\geq 0.\]
 Hence we deduce that
$f>0$ on the interval $I$ so that Theorem \ref{t2log} is immediate.

As mentioned above, the following recurrence relation will be needed
in the proof of Theorem \ref{t2log}.

\begin{lemma}\label{rr}
For $1\leq k\leq n$, we have
\begin{equation}\label{rec4}
d_n^{k-1}=(k+1)(n-k)d_n^{k+1}-(n-2k+1)d_n^k.
\end{equation}
\end{lemma}

\proof First, it is easy to
 establish the following recurrence relation for $1\leq k\leq n$,
\begin{equation}\label{rec5}
d_n^{k-1}=kd_n^k-d_{n-1}^{k-1}.
\end{equation}
By (\ref{rec1}) and (\ref{rec3}), we have
\begin{eqnarray*}
d_n^{k-1}&=&d_{n+1}^k-(n-k+1)d_n^k\\[5pt]
&=&(n+1)d_n^k-d_{n-1}^{k-1}-(n-k+1)d_n^k\\[5pt]
&=&kd_n^k-d_{n-1}^{k-1},
\end{eqnarray*}
as claimed. By (\ref{rec3}), (\ref{rec5}), for $1\leq
k\leq n$,  we find
\begin{eqnarray*}
d_n^k&=&(k+1)d_n^{k+1}-d_{n-1}^k\\[5pt]
&=&(k+1)d_n^{k+1}-\left(\frac{1}{n-k}d_n^k-\frac{1}{n-k}d_{n-1}^{k-1}\right)\\[5pt]
&=&(k+1)d_n^{k+1}-\frac{1}{n-k}d_n^k+\frac{1}{n-k}\left(kd_n^k-d_n^{k-1}\right)\\[5pt]
&=&(k+1)d_n^{k+1}-\frac{k-1}{n-k}d_n^k-\frac{1}{n-k}d_n^{k-1}.
\end{eqnarray*}
Consequently,
\[
d_n^{k-1}=(k+1)(n-k)d_n^{k+1}-(n-2k+1)d_n^k,
\]
as desired. \qed

In order to prove (\ref{bound}), we first give a lower bound
for $d_{n+1}^k/d_n^k$.

\begin{lemma}\label{lower} For $n\geq 1$ and $1\leq k\leq n-1$, we have
\begin{equation}\label{ne}
\frac{d_{n+1}^k}{d_n^k}\geq n+\frac{n-k}{n}.
\end{equation}
\end{lemma}

\proof  We proceed by induction on $n$. It is clear that
\eqref{ne} holds for $n=1$ and $n=2$. We now assume that \eqref{ne}
holds for positive integers less than $n$. By the recurrence
\eqref{rec1}, we have
\begin{align*}
\frac{d_{n+1}^k}{d_n^k}&=\frac{nd_n^k+(n-k)d_{n-1}^k}{d_n^k}\\[5pt]
&=n+(n-k)\frac{d_{n-1}^k}{d_n^k}\\[5pt]
&=n+(n-k)\frac{d_{n-1}^k}{(n-1)d_{n-1}^k+(n-k-1)d_{n-2}^k}.
\end{align*}
Thus \eqref{ne} can be recast as
\[
(n-1)+(n-k-1)\frac{d_{n-2}^k}{d_{n-1}^k}\leq n.
\]
So it  suffices to check that
\[
\frac{d_{n-1}^k}{d_{n-2}^k}\geq n-k-1.
\]
Since $n\geq 3$, by the inductive hypothesis,  we have
\begin{eqnarray*}
\frac{d_{n-1}^k}{d_{n-2}^k}&\geq & n-2+\frac{n-2-k}{n-2}\\[5pt]
&=&n-1-\frac{k}{n-2}\\[5pt]
&\geq &n-k-1.
\end{eqnarray*}
as required. \qed

Next we  give an upper bound for $d_{n+1}^k/d_n^k$.

\begin{lemma}\label{lup}
For $n\geq 4$ and $2\leq k\leq n-2$, we have
\begin{equation}\label{up}
\frac{d_{n+1}^k}{d_n^k}\leq n+\frac{n-k}{n}+\frac{n-k}{n^2}.
\end{equation}
\end{lemma}
\proof It follows from the recurrence \eqref{rec1} that
\begin{align*}
\frac{d_{n+1}^k}{d_n^k}&=n+(n-k)\frac{d_{n-1}^k}{d_n^k}\\[5pt]
&=n+(n-k)\frac{d_{n-1}^k}{(n-1)d_{n-1}^k+(n-k-1)d_{n-2}^k}.
\end{align*}
Thus \eqref{up} can be rewritten as
$$(n-1)+(n-k-1)\frac{d_{n-2}^k}{d_{n-1}^k}\geq \frac{n^2}{n+1},$$
that is,
\begin{equation}\label{r1}\frac{d_{n-1}^k}{d_{n-2}^k}\leq (n+1)(n-k-1).
\end{equation}
By recurrence \eqref{rec2} for $2\leq k\leq n-2$, we see that $$\frac{d_{n-1}^k}{d_{n-2}^k}\leq
n-1,$$
which implies (\ref{r1}). This completes the proof.  \qed

We are now ready to give the proof of Theorem \ref{t2log}.

\noindent{\it Proof.} It is easy to check that
 that the theorem holds for
$n=4,5,6$ and $2\leq k\leq n-2$. So we may assume that
$n\geq 7$.

We claim that the left hand side of (\ref{2log}) can be expressed as a cubic
function $f$ on $\frac{d_{n+1}^k}{d_n^k}$. By the recurrences
(\ref{rec1}), (\ref{rec2}), (\ref{rec3}) and (\ref{rec4}), we can
derive the following relations,
\begin{align*}
d_n^{k-2}&=(n-k+1)(n-k+3)d_n^k-(n-2k+3)d_{n+1}^k,\\[5pt]
d_n^{k-1}&=d_{n+1}^k-(n-k+1)d_n^k,\\[5pt]
d_n^{k+1}&=\frac{1}{(k+1)(n-k)}\left(d_{n+1}^k-kd_n^k\right),\\[5pt]
d_n^{k+2}&=\frac{1}{(k+1)(k+2)(n-k-1)(n-k)}\left((n-2k-1)d_{n+1}^k+(n+k^2)d_n^k\right).
\end{align*}
It follows that  \eqref{2log} can be rewritten as
$$A\cdot\left(C_3(n,k)\left(d_{n+1}^k\right)^3+C_2(n,k)\left(d_{n+1}^k\right)^2\left(d_n^k\right)
+C_1(n,k)\left(d_{n+1}^k\right)\left(d_n^k\right)^2+C_0(n,k)\left(d_n^k\right)^3\right)\geq
0,$$ where
\begin{align*}
A&=\frac{d_n^k}{(k+1)^2(n-k)^2(k+2)(n-k-1)},\\[5pt]
C_3(n,k)&=-n^2-5n+6k+6,\\[5pt]
C_2(n,k)&=n^3+n^2k+5n^2+3nk-10k^2+n-16k-6,\\[5pt]
C_1(n,k)&=n^2-2n+14k+14k^2+n^3+10nk^2-10n^2k-n^3k-3nk,\\[5pt]
C_0(n,k)&=-4n^2-12k^2-12k^3+10nk+18nk^2-9n^2k+n^2k^2-n^3k.
\end{align*}
Since $d_n^k$ are positive integers, it suffices to show that
\begin{equation}\label{c321}
C_3(n,k)\left(\frac{d_{n+1}^k}{d_n^k}\right)^3+C_2(n,k)\left(\frac{d_{n+1}^k}{d_n^k}\right)^2
+C_1(n,k)\left(\frac{d_{n+1}^k}{d_n^k}\right)+C_0(n,k)\geq 0.
\end{equation}
We now consider the function
$$f(x)=C_3(n,k)x^3+C_2(n,k)x^2+C_1(n,k)x+C_0(n,k),$$
with
\begin{equation}\label{diff}
f'(x)=3C_3(n,k)x^2+2C_2(n,k)x+C_1(n,k).
\end{equation}
We are going to show that $f'(x)<0$, for $2\leq x\leq n-1$.
As will be seen, the
quadratic function $f'(x)$ has a zero in the interval $[-1,k]$ and a zero in the
interval $[k,n]$.  At the point  $x=1$, we have \[f'(-1)=-(k+1)(n^3+12n^2-10nk+19n-34k-30).\]
Since for $n\geq 7$ and $2\leq k\leq n-2$, we find
\begin{align*}
\lefteqn{ n^3+12n^2-10nk+19n-34k-30}\\[5pt]
&\quad \geq n^3+12n(k+2)+19n-30-10nk-34k\\[5pt]
&\quad \geq (n^3-30)+2nk+(43n-34k)>0.
\end{align*}
This yields that $f'(-1)<0$.
Similarly, for $x=k$, we obtain that
\[
f'(k)=(k+1)(n-k)(n^2+n+2k-2)>0.
\]
Moreover, for $x=n$, we have
\begin{equation}\label{diff2}
f'(n)=-(n-k)(n^3+4n^2-10nk+14k-21n+14).
\end{equation}
To prove $f'(n)<0$, it is sufficient to show that for $2\leq k\leq n-2$,
\[
n^3+4n^2-10nk+14k-21n+14>0.
\]
We have two cases for the ranges of $k$.
For $2\leq k\leq n-3$, we have
$$n^3+4n^2-10nk+14k-21n+14=n\left((n-3)^2+10(n-k-3)\right)+14k+14>0,$$
Meanwhile, for $k=n-2$,
$$n^3+4n^2-10nk+14k-21n+14=n(n-3)^2+4n-14>0.$$ Thus $f'(n)<0$ is valid
for $2\leq k\leq n-2$.
Then we reach the conclusion that
 $f'(x)$ has a zero in the interval $[-1,k]$ and
 a zero in the interval$[k,n]$.

We continue to demonstrate that $f'(x)<0$ in the interval $I$.
By Lemma \ref{ne}, for $k\geq
2$ we have  $$\frac{d_{n+1}^k}{d_n^k}\geq n+\frac{n-k}{n}>n,$$ which means  that $f'(x)$ has no zero on the interval $I$. Since $n\geq
k+2$, it is easily seen that
\begin{align*}
C_3(n,k)&=-(n^2+5n-6k-6)\\[5pt]
& \leq -\left((k+2)^2+5(k+2)-6k-6\right)\\[5pt]
& \leq -(k^2+3k+8)<0.
\end{align*}
Since $f'(n)<0$, we see that $f'(x)<0$ in the interval $I$, as expected. In other words, $f(x)$ is strictly
decreasing on this interval.

Up to now, we have shown that $f(x)$ is strictly decreasing on the
interval $I=[n+\frac{n-k}{n},n+\frac{n-k}{n}+\frac{n-k}{n^2}]$.
So it remains to  prove that
$$f\left(n+\frac{n-k}{n}+\frac{n-k}{n^2}\right)>0.$$
Since
$$f\left(n+\frac{n-k}{n}+\frac{n-k}{n^2}\right)=\frac{h(k)(n-k)^2}{n^6},$$
where
\begin{align*}
h(k)&=(-10n^4-26n^3-28n^2-18n-6)k^2+(-n^6+20n^5+27n^4+19n^3-7n-6)k\\[5pt]
\quad\quad & \quad \quad +(n^7-10n^6-4n^5-4n^4+9n^3+7n^2+6n).
\end{align*}
Clearly, the proof will be complete as long as we can show that $h(k)\geq
0$ for $n\geq 7$ and $2\leq k\leq n-2$.

 Regard $h(x)$ as a continuous function on $x$, that is,
 \begin{align*}
h(x)&=(-10n^4-26n^3-28n^2-18n-6)x^2+(-n^6+20n^5+27n^4+19n^3-7n-6)x\\[5pt]
\quad\quad &\quad \quad +(n^7-10n^6-4n^5-4n^4+9n^3+7n^2+6n).
\end{align*}
Since the leading coefficient $-10n^4-26n^3-28n^2-18n-6$ of $h(x)$ is negative,
we only need to prove that $h(2)>0$ and
$h(n-1)>0$. For $n\geq 7$, we have
\begin{align*}
h(n-1)&=n(n^5-3n^4+2n^3+2n^2+2n+1)\\[5pt]
&=n\left(n^3(n-1)(n-2)+2n^2+2n+1\right)>0,
\end{align*}
and
\begin{align*}
h(2)&=n^7-12n^6+36n^5+10n^4-57n^3-105n^2-80n-36\\[5pt]
&=n^5(n-5)(n-7)+n^4(n-6)+16n^3(n-7)+55n^2(n-7)\\[5pt]
\quad \quad &\quad \quad +80n(n-1)+200n^2-36>0.
\end{align*}
In summary, we have confirmed  that $h(k)>0$ for $n\geq 7$ and $2\leq k\leq
n-2$. This completes the proof.\qed

\section{The reverse ultra log-concavity}

This section is concerned with the reverse ultra log-concavity of
$d_n^k$. Recall that sequence $ \{a_k\}_{0\leq k\leq n}$ is called
ultra log-concave if $\left\{a_k\big/{n\choose k}\right\}$ is
log-concave; see Liggett \cite{Lig}.  This condition can be restated
as
\begin{equation}\label{ultralc}
k(n-k)a_k^2-(n-k+1)(k+1)a_{k-1}a_{k+1}\geq 0.
\end{equation}
It is well known that if a polynomial has only real zeros, then its
coefficients form an ultra log-concave sequence. As noticed by
Liggett \cite{Lig},  if a sequence $\{a_k\}_{0\leq k\leq n}$ is
ultra log-concave, then the sequence $\{k!a_k\}_{0\leq k\leq n}$ is
log-concave.

In comparison with ultra log-concavity, a
 sequence is said to be reverse ultra log-concave if it satisfies
the reverse relation of \eqref{ultralc}, that is,
\begin{equation}\label{reverse}
k(n-k)a_k^2-(n-k+1)(k+1)a_{k-1}a_{k+1}\leq 0.
\end{equation}
Chen and Gu \cite{Chen} have shown the Boros-Moll polynomials have
this reverse ultra log-concave property. We shall show that the
sequence $\{d_n^k\}_{0\leq k\leq n}$ also possesses this property.

\begin{theorem}
For $1\leq k\leq n-1$, we have
\[
\frac{d_n^{k-1}}{{n\choose k-1}}\cdot\frac{d_n^{k+1}}{{n\choose
k+1}}\geq \left(\frac{d_n^k}{{n\choose k}}\right)^2,
\]
or equivalently,
\begin{equation}\label{reu}
(n-k+1)(k+1)d_n^{k-1}d_n^{k+1}\geq k(n-k)\left(d_n^k\right)^2.
\end{equation}
\end{theorem}

\proof According to the recurrence  relations \eqref{rec3} and \eqref{rec4},
we find that \eqref{reu} can be reformulated as
\begin{equation}\label{eqq}
(n-k+1)\left(\frac{d_{n+1}^k}{d_n^k}\right)^2-(n-k+1)(n+1)\left(\frac{d_{n+1}^k}{d_n^k}\right)+k(2n-2k+1)\geq
0.
\end{equation}
The discriminant of the quadratic polynomial of the left side of
\eqref{eqq} in $d_{n+1}^k/d_n^k$ equals
$$\Delta=\left((n-k+1)(n+1)\right)^2-4k(n-k+1)(2n-2k+1).$$
We claim that $\Delta>0$ for $1\leq k\leq n-1$. Put
\[ f(k)=\Delta=8k^2-(n^2+10n+5)k+(n^3+3n^2+3n+1).\]
 Since $n\geq
k+1$, we have
\begin{align*}
f'(k)=&16k-(n^2+10n+5)\\[5pt]
=&-(n^2+10n-16k+5)\\[5pt]
\leq & -\left((k+1)^2+10(k+1)-16k+5\right)\\[5pt]
=&-(k-2)^2-12<0,
\end{align*}
which implies that $f(k)$ is  monotone decreasing for $1\leq k\leq
n-1$. Furthermore, \[ f(n-1)=2\left((n-2)^2+3\right)>0.\]¡¡Thus,
$\Delta>0$  for  $ 1\leq k\leq n-1$.
Consequently, the quadratic function has two distinct real zeros. If
we can show that for $1\leq k\leq n-1$, $d_{n+1}^k/d_n^k$ is larger than the maximal zero, then (\ref{eqq})
holds since $n-k+1>0$. Thus we still have to  show that
\begin{equation}\label{eq2}
\frac{d_{n+1}^k}{d_n^k}>
\frac{(n-k+1)(n+1)+\sqrt{\Delta}}{2(n-k+1)}
=\frac{n+1}{2}+\frac{\sqrt{\Delta}}{2(n-k+1)}
\end{equation}
In view of \eqref{ne}, we see that \eqref{eq2} can be deduced
 from the following inequality
$$n+\frac{n-k}{n}\geq \frac{n+1}{2}+\frac{\sqrt{\Delta}}{2(n-k+1)},$$
which is equivalent to
$$(n-k+1)(n^2+n-2k)\geq n\sqrt{\Delta}.$$
Since both sides are positive, we can transform the above relation into the following form
$$\left((n-k+1)(n^2+n-2k)\right)^2\geq n^2\Delta.$$
Evidently,
\begin{align*}
&\ \left((n-k+1)(n^2+n-2k)\right)^2-n^2\Delta\\[5pt]
&=(n-k+1)\left(4n^2k(2n-2k+1)-4k(n-k+1)(n^2+n-k)\right)\\[5pt]
&=4k(n-k+1)(n-k)(n^2-n+k-1)\geq 0,
\end{align*}
for $1\leq k\leq n-1$. This completes the proof.  \qed

\vspace{.2cm} \noindent{\bf Acknowledgments.} This work was
supported by  the 973 Project, the PCSIRT Project of the Ministry of
Education,  and the National
Science Foundation of China.

\end{document}